\documentclass[11pt,a4paper]{article}
\usepackage{graphicx}
\usepackage{amsfonts}
\usepackage{latexsym}
\usepackage{epsfig}

\newcommand{\beql}[1]{\begin{equation}\label{#1}}
\newcommand{\eeq}{\end{equation}}
\newcommand{\comment}[1]{}

\newcommand{\eqref}[1]{{\rm (\ref{#1})}}

\newcommand{\Abs}[1]{{\left|{#1}\right|}}

\newcommand{\Norm}[1]{{\left\|{#1}\right\|}}

\newcommand{\Ceil}[1]{{\left\lceil{#1}\right\rceil}}
\newcommand{\Qed}{\ \\\mbox{$\Box$}}

\newcommand{\Set}[1]{{\left\{{#1}\right\}}}
\newcommand{\given}{{\ |\ }}

\newcommand{\RR}{{\mathbb R}}

\newcommand{\ZZ}{{\mathbb Z}}
\newcommand{\NN}{{\mathbb N}}
\newcommand{\TT}{{\mathbb T}}
\newcommand{\QQ}{{\mathbb Q}}

\newcommand{\inner}[2]{{\langle #1, #2 \rangle}}

\newcommand{\supp}{{\rm supp\,}}

\newcommand{\ft}[1]{\widehat{#1}}

\newcommand{\Zset}{{\cal Z}}
\newcommand{\MM}{{\cal M}}

\newcommand{\Msoz}{{\cal M}^*(\Omega,z)}

\newcounter{open}
\newcounter{dfn}
\def\thedfn{\arabic{dfn}}

\newcounter{obs}
\def\theobs{\arabic{obs}}

\newcounter{thm}
\def\thethm{\arabic{thm}}

\newcounter{othm}
\def\theothm{\Alph{othm}}

\newcounter{mysec}
\newcounter{mysubsec}[mysec]
\newtheorem{theorem}{Theorem}
\newtheorem{corollary}{Corollary}
\newtheorem{lemma}{Lemma}
\newtheorem{proposition}{Proposition}

\newtheorem{remark}{Remark}

\newtheorem{problem}{Problem}

\newcounter{rem}
\newcounter{rev}
\begin{document}

\title{On pointwise estimates of positive definite functions with given support}
\author{Mihail N. Kolountzakis\footnote{
Supported in part by European Commission IHP Network HARP
(Harmonic Analysis and Related Problems),
Contract Number: HPRN-CT-2001-00273 - HARP.} \and
Szil\'ard Gy. R\'ev\'esz
\thanks {The second author was supported in part by the
Hungarian National Foundation for Scientific Research,
Grant \# T034531 and T032872.}}

\date{January 2003}

\maketitle

\begin{abstract}
The following problem originated from a question due to Paul
Tur\' an. Suppose $\Omega$ is a convex body in Euclidean space
$\RR^d$ or in $\TT^d$, which is symmetric about the origin. Over
all positive definite functions supported in $\Omega$, and with
normalized value $1$ at the origin, what is the largest possible
value of their integral? From this Arestov, Berdysheva and Berens
arrived to pose the analogous pointwise extremal problem for
intervals in $\RR$. That is, under the same conditions and
normalizations, and for any particular point $z\in\Omega$, the
supremum of possible function values at $z$ is to be found.
However, it turns out that the problem for the real line has
already been solved by Boas and Kac, who gave several proofs and
also mentioned possible extensions to $\RR^d$ and non-convex
domains as well.

We present another approach to the problem, giving the
solution in $\RR^d$ and for several cases in $\TT^d$.
In fact, we elaborate on the fact that the problem is
essentially one-dimensional, and investigate non-convex open domains as well.
We show that the extremal problems are equivalent to more familiar ones over
trigonometric polynomials, and thus find the extremal values for a few
cases. An analysis of the relation of the problem for the space
$\RR^d$ to that for the torus $\TT^d$ is given, showing that the former case is
just the limiting case of the latter. Thus the hiearachy of difficulty is
established,
so that trigonometric polynomial extremal problems gain recognition again.
\end{abstract}

\vskip 1cm



\vfill
\eject

{\bf \S 1 Extremal Problems for positive definite functions, periodic and not}

\vskip0.6cm

Let us denote $\TT^d :=\left[-{1\over 2},{1 \over 2}\right)^d\subset \RR^d$
with the usual
modified topology of periodicity, that is, take the topology of
$\TT^d:=\RR^d/\ZZ^d$. For $\Omega\subseteq\TT^d$ any open domain, we put

\beql{Fstar}
{\cal F}^* (\Omega):=\{f:\TT^d\to\RR\,:\, \supp f \subseteq \Omega, f(0)=1, f\,
{\rm positive}\,\,{\rm definite} \},
\eeq
and, analogously, when $\Omega\subseteq\RR^d$ is any open set,

\beql{Fclass}
{\cal F}(\Omega):=\{f:\RR^d\to\RR\,:\, \supp f \subseteq \Omega, f(0)=1, f\,
{\rm positive}\,\,{\rm definite} \}.
\eeq

Recall that positive definiteness of functions (and even measures
and tempered distributions) can be defined or equivalently
characterized by nonnegativity of Fourier transform. In case of
\eqref{Fstar} positive definiteness means $\ft{f}(n)\ge 0 \,
(\forall n \in \ZZ^d ) $, while in case of \eqref{Fclass} it
means $\ft{f}(x)\ge 0 \, (\forall x \in \RR^d ) $.

In 1970 in a discussion with S. B. Stechkin \cite{Stechkin} Paul Tur\' an
posed the following problem. Let $d=1$ and $\Omega:=(-h,h)\subset\TT$.
What is the largest possible value of the integral $\int_{\TT} f$ over all
$f\in {\cal F}^*((-h,h))$? The question was later investigated in higher
dimensions and in $\RR^d$ as well. As a natural condition for the above
Tur\'an extremal problem, convexity of the underlying domain $\Omega$ is
usually supposed.

For an account of the problem see the papers \cite{AB} and
\cite{KR} and the references therein. However, no authors seem to
have noticed that already Boas and Kac settled the analogous (and
relatively easy) case of an interval $(-h,h)\subset\RR$, see
Theorem 5 of \cite{BK}.

In \cite{ABB} the natural pointwise analogue of the above question of Tur\'an
was studied for intervals in $\TT$ or $\RR$. For general domains in arbitrary
dimension these problems can be formulated as follows.

\begin{problem}\label{pr:pointwise}
{\rm (Boas-Kac - type pointwise extremal problem for the space)}.
Let $\Omega \subseteq \RR^d$ be an open set,
and let $f:\RR^d\to\RR$ be a
positive definite function with $\supp f \subseteq \Omega$ and
$f(0)=1$. Let also $z\in\Omega$. What is the largest possible
value of $f(z)$? In other words, determine

\beql{omegaz} {\cal M}(\Omega, z) := \sup_{f\in {\cal F}(\Omega)}
f(z). \eeq
\end{problem}

\begin{remark}Obviously, ${\cal M}(\Omega, z) \le 1$, as
$1 \pm  f(z) =\int_{\RR} (1\pm \exp(2\pi i zt)) \ft{f}(t) dt$
$=\int_{\RR} (1\pm \cos(2\pi zt)) \ft{f}(t) dt \ge 0.$
\end{remark}

One might miss a more precise specification of the function class
$f:\RR^d\to\RR$ here and similarly in the problems listed below.
The fact that considering $L^1$, $C$ or $C^\infty$ leads to the same
answer ie. same extremal values, will be discussed at the beginning
of \S 2.

\begin{problem}\label{pr:torus}
{\rm (Tur\'an - type pointwise extremal problem for the torus)}.
Let $\Omega \subseteq \TT^d$ be any open set,
and let $f:\TT^d\to\RR$ be a
positive definite function with $\supp f \subseteq \Omega$ and
$f(0)=1$. Let also $z\in\Omega$. What is the largest possible
value of $f(z)$? In other words, determine

\beql{Mstar} {\cal M}^* (\Omega, z) := \sup_{f\in {\cal
F}^*(\Omega)} f(z). \eeq

\end{problem}

\begin{remark}
Let $\Omega \subseteq (-{1\over 2}, {1\over 2})^d$ and
$f:\Omega\to\RR$.
For the function $f$ to be positive definite {\em on the
torus} means a nonnegativity condition for the Fourier Transform
$$
\ft{f}(\xi) = \int_{\RR^d} e^{2\pi i \inner{\xi}{x}} f(x)~dx
$$
only for a discrete set of values of $\xi$, namely $\xi \in \ZZ^d$,
while positive definiteness of $f$ as a function
on $\RR^d$ is equivalent to nonnegativity of the
Fourier transform $\widehat f$ for all occurring values.
From this it follows that we always have
\begin{equation}\label{comparison}
{\cal M}^*(\Omega, z) \ge {\cal M}(\Omega, z).
\end{equation}
\end{remark}

The extremal value in the above Problem \ref{pr:pointwise} was
estimated together with its periodic analogue Problem
\ref{pr:torus} in the work \cite{ABB} for dimension $d=1$.
However, already Boas and Kac has solved the $d=1$ case of
Problem \ref{pr:pointwise}, which seem to have been unnoticed in
\cite{ABB}.

These problems are not only analogous, but also related to each other, and, in
fact, Problem \ref{pr:pointwise} is only a special, limiting case of the more
complex Problem \ref{pr:torus} (see Theorem \ref{limitcase} below).
On the other hand, already Boas and Kac observed, that Problem \ref{pr:pointwise}
(dealt with for $\RR$ in \cite{BK}) is connected to trigonometric polynomial
extremal problems. In particular, from the solution to the interval case they
deduced the value \eqref{cosvalue} below of the extremal problem due to Carath\'eodory
\cite{Cara} and Fej\'er \cite{Fej}.
They also established a connection (see \cite[Theorem 6]{BK}) what
corresponds to the one-dimensional case of the first part of our Theorem
\ref{th:connection}.

It is appropriate at this point to consider also the following type
of trigonometric polynomial extremal problems. Denote for any
$H\subseteq \NN_2:=\NN\cap \left[2,\infty\right)$

\beql{Ficlass}
\Phi(H) := \Set{\varphi:\TT\to\RR_+\,:\, \lambda\in\RR, \varphi\ge 0, \varphi (t) \sim 1 + \lambda
    \cos 2\pi t + \sum_{k\in H} c_k \cos 2\pi kt }
\eeq
and with a given $m\in\NN_2$ and $H\subseteq\NN_2$ also

\beql{Fimstar}
\Phi_m (H) := \{\varphi:\TT\to\RR \, : \, \, \lambda\in\RR, \varphi (\frac{j}{m})
\ge 0 \, (j\in \ZZ), \, \, \varphi (t) = \qquad \qquad
\qquad \qquad
\eeq
$$
\qquad \qquad \qquad \qquad \qquad \qquad \qquad \qquad = 1 +
\lambda \cos 2\pi t + \sum_{k\in H} c_k \cos 2\pi kt \} .
$$

\begin{problem}\label{pr:trigo}
{\rm (Carath\'eodory-Fej\'er type trigonometric polynomial problem)}.
Determine the extremal quantity
\beql{Mextr}
M(H):= \sup \{\lambda= 2 \ft \varphi (1)\,:\, \varphi \in \Phi(H)\}.
\eeq
\end{problem}

\begin{remark} Observe that $M(H)\le 2$, always, as
$$
\Abs{\lambda/2} = \Abs{\ft\varphi(1)} \le \Norm{\varphi}_1 =
 \int \varphi = \ft\varphi(0) = 1.
$$
\end{remark}

\begin{problem}\label{pr:fft}
{\rm (Discretized Carath\'eodory-Fej\'er type extremal problem)}.
Determine \beql{Mmstar} M_m(H):= \sup \{\lambda= 2 \ft \varphi
(1)\,:\, \varphi \in \Phi_m(H)\}. \eeq

\end{problem}

\begin{remark} It should be remarked here that obviously we have
$\Phi(H)\subseteq\Phi_m(H)$.
So we always have $M_m(H) \ge M(H)$.
\end{remark}

In this note we present the exact solution of Problem \ref{pr:pointwise} in
line of what the paper \cite{BK} suggests. In fact, we have to acknowledge
that Boas and Kac mentioned the possibility of extending one of their methods
-- Poisson summation -- to higher dimensions, so some parts of what follows can
be interpreted as implicitly present already in their work \cite{BK}. But
here we obtain some results also for the more complex periodic version.

However, the main result of the present investigation is perhaps the
understanding that the above point-value extremal problems are in fact equivalent
to the above trigonometric polynomial extremal problems, thus transferring
information on one problem to the equivalent other problem in several cases.
Until now the equivalence formulated below remained unclear in spite of the
fact that, e.g., Boas and Kac found ways to deduce the solution of the
trigonometric extremal problems in  Problem
\ref{pr:trigo} from their results on Problem \ref{pr:pointwise}.
We also obtain a clear picture of the limiting relation between torus problems
and space problems, and, parallel to this, between the finitely conditioned
trigonometric polynomial extremal problems of Problem \ref{pr:fft} and the
positive definite trigonometric polynomial extremal problems of Problem
\ref{pr:trigo}.

\vskip0.6cm

{\bf \S 2 Preliminaries. Formulation of the Equivalence Results.}

\vskip0.6cm

Note that in the above definitions  \eqref{Fstar}, \eqref{Fclass}
or \eqref{Ficlass}, \eqref{Fimstar} it is left a bit unclear,
what function classes are considered as $\RR^d\to\RR$,
$\TT^d\to\RR$ or $\TT\to\RR$. However, this causes no ambiguity,
since it is not hard to see that the extremal problems
\eqref{omegaz}, \eqref{Mstar}, \eqref{Mextr} or \eqref{Mmstar}
yield the same extremal values when e.g., integrable functions
(with continuity of $f$ supposed only at $z$ in case of
\eqref{omegaz} or \eqref{Mstar}) are considered, and when e.g.,
compactly supported $C^\infty$ functions are taken into account.
Indeed, on $\TT$ or $\TT^d$ this follows after a convolution by
the Fej\'er kernels, say. The same way we can restrict ourselves
even to trigonometric polynomials in $\Phi(H)$ or $\Phi_m(H)$ as
well.

Passing on to the case of the real space $\RR^d$, first we show
that it suffices to consider {\it bounded} open sets only. To
this end let us consider the auxiliary positive definite function

\beql{triangle}
\Delta_R(x):= {1 \over |B_{R/2}|} \chi_{B_{R/2}}* \chi_{B_{R/2}}
\eeq
with $B_r:=\{x\in \RR^d \,: \, |x|\le r\}$, and take $
f_N:=f\Delta_N$ to obtain

$$
{\cal M}(\Omega,z)=\lim_{N\to \infty}{\cal M}(\Omega_N,z)
=\lim_{N\to \infty}{\cal M}({\rm int}\, \Omega_N,z),
$$
where $\Omega_N:=\{x\in\Omega\, :\, |x|\le N\}=\Omega\cap B_N$,
and thus $\Omega_N\subseteq {\rm int}\, \Omega_{N+1}$.

Next observe that for any bounded open $\Omega$, the condition
$\supp f \subseteq \Omega$ entails that $\supp f$ is compact and of
a fixed positive distance $\eta$ from the boundary of $\Omega$.
Thus convolution of $f$ with the (convolution) square of some
approximate identity $k_\delta$ with $\supp k_\delta\subseteq
B_\delta$ leads to a function $f_\delta:=f * k_\delta * k_\delta$
satisfying $\supp f_\delta \subseteq \supp f + B_{2\delta} \subseteq
\Omega$ if $\delta < {1\over 2}\eta$. Hence with a smooth
$k_\delta$ we have $f_{\delta}\in {\cal F}(\Omega)\cap
C^\infty(\Omega)$, while for arbitrary fixed $\epsilon
> 0$ and with $\delta$ correspondingly small enough $f_\delta
(z)\ge f(z) - \epsilon$ in view of the continuity of $f$ at $z$.

Now let us define for $z\in\Omega$ the derived set

\beql{Homegaz}
H(\Omega,z):=\{k\in\NN_2\, : \, kz\in\Omega,\, -kz\in\Omega\}
\eeq

Our first goal is to show that in fact the Boas-Kac type
Problem \ref{pr:pointwise}
is a one-dimensional problem. This is contained in the following result.

\begin{theorem}
\label{th:connection}
Let $0\in\Omega\subseteq\RR^d$ be any open set and
$z \in \Omega \cap (-\Omega)$. With the above notations we have

$$
{\cal M} (\Omega,z)=\frac12 M(H(\Omega,z)).
$$

\end{theorem}

\begin{remark} Note that in case $z\in\Omega$, $z\notin-\Omega$,
we trivially conclude ${\cal M}(\Omega,z)=0$ since for all
$f\in{\cal F}(\Omega)$, $\supp f \subseteq \Omega\cap (-\Omega)$
follows from \eqref{symmetry} below. Also $0\in\Omega$ is
necessary, for a positive definite function $f$ must vanish a.e.
if $0 \notin \supp f$.
\end{remark}

To tackle the Tur\'an-type Problem \ref{pr:torus}, one may consider
$f\in L^1(\TT^d)$ with continuity supposed at $z$, or even
$f\in C^\infty(\TT^d)$.

Here positive definiteness of $f$ is equivalent to $\ft{f}(n)\ge
0 \quad (\forall n \in \ZZ^d)$, and similarly to
\eqref{symmetry}, one gets $f(x)=f(-x) \quad (\forall
x\in\TT^d)$. Thus $\supp f$ is symmetric, hence $\supp f \subseteq
\Omega\cap (-\Omega)$.

Once again we see that \eqref{Mstar}
vanishes unless $z\in \Omega\cap (-\Omega)$ and that it suffices
to restrict ourselves to sets symmetric about the origin. In
other words, if $z\notin \Omega$ or if $z\notin (-\Omega)$, then
${\cal M}^*(\Omega, z)=0$, while for $z=0$ obviously ${\cal
M}^*(\Omega, 0)=1$. These are the trivial cases, and for the
remaining cases we introduce a further notation. Put
\beql{Zset} {\cal Z}:= {\cal Z}(z):= \{nz \,({\rm mod}\,\TT^d)
\given n\in\ZZ\}. \eeq

The set $\cal Z$ is finite if and only if we have $z\in\QQ^d$,
that is, $z=({p_1 \over q_1},\dots, {p_d \over q_d})$ with $p_j
,\, q_j \in \ZZ, \,(p_j,q_j)=1 \,\, (j=1,\dots,d)$. In this case
we have with $m=[q_1,\dots,q_d]$, the least common multiple of the
denominators, that $mz=0 \, ({\rm mod}\,\TT^d)$, and for arbitrary
$n, n' \in \ZZ$ $nz=n'z \, ({\rm mod}\,\TT^d)$ if and only if
$n\equiv n' \,({\rm mod}\, m)$.

Let us keep the definition \eqref{Homegaz} with an interpretation
$\pmod {\TT^d}$ for infinite $\cal Z$. On the other hand, in case
\#${\cal Z}=m$ we put

\beql{Hmstar} H_m(\Omega,z):=\{k\in [2,m/2]\, : \, kz\in\Omega,\,
-kz\in\Omega\} =H(\Omega,z)\cap [2,m/2]. \eeq
Moreover, for any set $H\subset\ZZ$ we define
$$
H(m):=\{k\in [2,m/2]\,:\, \exists h\in H \mbox{such that} \pm k
\equiv h \pmod m \}.
$$

\begin{remark}
Note the following relations for an arbitrary $H\subseteq\NN_2$. 
First, if there exists any index $k\in H$ with $k
\equiv 1 \pmod m$, then we obtain $M_m(H)=\infty$, since $1+a\cos
{2\pi t} - a\cos {2 k \pi t}$ is nonnegative at $j/m$ for all
$j=1,\dots,m$ and any $a \in \RR$. In fact, for $k\equiv \ell
\pmod m$ obviously $\cos {2 k \pi t} - \cos {2 \ell \pi t}$ is
vanishing at all points of the form $j/m$, hence the coefficients
can be changed mod $m$ to reduce $\varphi$ to a trigonometric
polynomial of degree at most $m$. Moreover, since this can be
used even for negative indices, and as $\cos{-k 2 \pi t}=\cos{ k
2 \pi t}$, in fact we can reduce the support of $\ft\varphi$ to
$[0,m/2]$. That is, either $M_m(H)=\infty$ (in case there is a
$k\in H$ with $k\equiv \pm 1 \pmod m$), or $M_m(H)=M_m(H(m))$.
\end{remark}

Now we can formulate

\begin{theorem}\label{th:torus}
Let $0\in\Omega\subseteq\TT^d$ be any open set and
$z \in \Omega \cap (-\Omega)$. Then the extremal quantity
\eqref{Mstar} depends only on the set ${\cal Z}$. In case ${\cal
Z}$ is infinite, we have
\beql{denseformula}
{\cal M}^* (\Omega,z)=\frac 12 M(H(\Omega,z)).
\eeq
In case $\#{\cal Z}=m$ is finite, we have

\beql{discreteformula} {\cal M}^* (\Omega, z) = \frac 12
M_m(H_m(\Omega,z)). \eeq
\end{theorem}

\vskip0.6cm

{\bf \S 3 Proof of Theorem \ref{th:connection}.}

\vskip0.6cm

First note that it suffices to consider symmetric sets
$\Omega'=\Omega\cap (-\Omega)$ only. Indeed, if $\Omega$ is arbitrary, and
$f\in{\cal F}(\Omega)$, $f\in C^\infty_0(\RR^d)$,
then by $\ft f \ge 0$ Fourier inversion yields

\beql{symmetry}
f(x)=\overline{f(x)}=\overline{\int \ft f(y)e^{2\pi i\inner{x}{y}}dy}
=\int \ft f (y)e^{-2\pi i\inner{x}{y}}dy=f(-x).
\eeq

Thus for all $f\in{\cal F}(\Omega) \quad \supp f$ is necessarily
symmetric. On the other hand, $H(\Omega,z)$ is symmetrized by
definition \eqref{Homegaz} with respect to $\Omega$. Hence we can
restrict ourselves to symmetric sets. Without loss of generality
we can assume that $\Omega$ is also bounded.

Now given a bounded symmetric open set $\Omega$
the proof consists of proving the two inequalities below.

\vskip0.4cm

\noindent
\underline{ ${\cal M}(\Omega, z) \le M(H(\Omega,z))/2$ }

\vskip0.4cm

Let $f$ have $f(0)=1$, be positive definite and have support in $\Omega$.
Define also the positive definite Radon measure
$$
\mu_z := \sum_{k\in\ZZ} \delta_{kz}.
$$
The function $f$ being continuous, the measure
\beql{nu-measure}
\nu_z = f\cdot \mu_z = \sum_{k\in\ZZ} f(kz) \delta_{kz}
\eeq
is well defined and positive definite as well.

Notice now, because of the boundedness of $\Omega$, that the sum
in \eqref{nu-measure} is actually a finite one. More precisely, if we have
e.g., $\Omega \subseteq B_n$, then we find

$$
\nu_z := \sum_{k=-(n-1)}^{n-1} f(kz) \delta_{kz}
=\delta_0+f(z)(\delta_z+\delta_{-z})+\sum_{k\in H(\Omega,z)}
f(kz)(\delta_{kz}+\delta_{-kz}) ,
$$
and that
$$
0 \le \widehat{\nu_z}(x) = 1 + 2f(z)\cos{2\pi\inner{z}{x}} +
  \sum_{k\in H(\Omega,z)} 2 f(kz) \cos{2\pi k \inner{z}{x}}
,\ \ \ (x \in \RR^d).
$$
Setting $t = \inner{z}{x}$ and observing that the trigonometric polynomial
$$
1 + 2f(z)\cos{2\pi t} + \sum_{k\in H(\Omega,z)} 2 f(kz) \cos{2\pi k t}
$$
is nonnegative, we obtain $2 f(z) \le M(H(\Omega,z))$.

\vskip0.4cm

\noindent
\underline{ ${\cal M}(\Omega, z) \ge M(H(\Omega,z))/2$}

\vskip0.4cm

For a function $\varphi\,:\,\TT\to\RR$ let us call the {\it
(restricted) spectrum} of $\varphi$ the set $S:=S(\varphi):=\supp
\ft{\varphi}\cap\NN_2\subseteq \NN_2$. Also, we will use the term
{\it full spectrum} and the notation $S':=S'(\varphi)$ for the set
$S':=\{-1,0,1\}\cup S \cup (-S)$, whether the exponential Fourier
coefficients at $-1,0$ or $1$ happen to vanish or not.

Take any trigonometric polynomial $\varphi \in \Phi (H)$ with
spectrum $S\subseteq H:=H(\Omega,z)$. Recall that taking the
supremum in \eqref{Mextr} over the function class \eqref{Ficlass}
yields the same result as considering such trigonometric
polynomials only.

Consider the measure
$$
\alpha_z := \delta_0 + (\lambda/2)(\delta_z + \delta_{-z}) +
\sum_{k\in S
} (c_k/2) (\delta_{kz} + \delta_{-kz}),
$$
whose Fourier transform is essentially equal to the polynomial
$\varphi (t)$ in \eqref{Ficlass}. Hence $\alpha_z$ is a positive
definite measure.

Take now the ``triangle function'' $\Delta_{\epsilon}$ defined as in
\eqref{triangle}, but here with a subscript $\epsilon$ small enough
to guarantee that
\begin{enumerate}
\item
The sets $kz + B_{\epsilon}$, $k \in S'$, are disjoint,
i.e., $\epsilon < {|z|\over 2}$, and
\item
These sets are all contained in $\Omega$,
i.e., $\epsilon < {\rm dist}\{\partial \Omega, S'z \}$.
\end{enumerate}
Finally define
$$
f := \alpha_z * \Delta_\epsilon ,
$$
which is a positive definite function supported in $\Omega$ with
value $1$ at the origin and with $f(z) = \lambda / 2$. This
proves that ${\cal M}(\Omega, z) \ge M(H(\Omega,z))/2$, as
desired. \Qed

\vskip0.6cm

{\bf \S 4 Applications of Theorem \ref{th:connection}}

\vskip0.6cm

The first application concerns the original convex case of the
pointwise Boas-Kac type problem formulated in Problem 1. A
symmetric, bounded convex domain with nonempty interior -- that
is, a {\it convex body} -- defines a norm. So for a vector $x$
let $||x||$ denote the norm of $x$ defined by $\Omega$, that is
$$
||x|| : = \inf\Set{\lambda>0:\ {1\over\lambda}x \in \Omega}.
$$
In other words, $\Omega$ is the unit ball of the norm $||\cdot||$.

\begin{corollary}
\label{th:answer-to-pointwise-problem} {\bf (Boas -- Kac
\cite{BK}).} Let $\Omega \subseteq \RR^d$ be a convex open
domain, symmetric about $0$. Suppose that
\beql{interval}
{1 \over n+1} \le ||z|| < {1 \over n},
\eeq
for some $n \ge 1$. Then
$$
{\cal M}(\Omega, z) = \cos{\pi\over n+2}.
$$
\end{corollary}

\noindent {\bf Proof of Corollary \ref{th:answer-to-pointwise-problem}.}
First observe that for the symmetric, convex, bounded, open set
$\Omega$ the norm of $z$ satisfies \eqref{interval} if and only
if $H(\Omega,z)=[2,n]$. Thus by Theorem \ref{th:connection} the
problem reduces to the extremal problem \beql{ntrigproblem}
M_n:=\sup \{\lambda \, :\, \exists \varphi (t) \ge 0, \varphi (t)
= 1 + \lambda \cos 2\pi t + \sum_{k=2}^{n} c_k \cos 2\pi kt \}.
\eeq This problem was settled by Fej\'er, see e.g., \cite{Fej} or
\cite[p. 869-870]{Fgesamm}. To finish the proof, we quote from
these or from \cite[Problem VI. 52, p.\ 79]{polya-szego} the
formula
\beql{cosvalue}
M_n = 2 \cos{\pi \over n+2}.
\eeq
\Qed

Note that \cite[Theorem 2]{ABB} gave the estimate $\frac{n}{n+1}\le
{\cal M}(\Omega,z) \le \frac12 (1+\cos(\frac{\pi}{n+1}))$
for the one-dimensional case. The above exact solution and some calculation
shows that both of these estimates are sharp for $n=1$, but none of them
is for $n>1$. However, this is covered (at least for $d=1$) by
\cite[Theorem 2]{BK}.

Now the $n\to\infty$ limiting case easily leads to

\begin{corollary}\label{BoasKac} {\bf (Boas -- Kac \cite{BK}).}
Suppose that the open set  $\Omega\subseteq \RR^d$ contains all
integer multiples of the point $z\in\RR^d$. Then $ {\cal
M}(\Omega,z)=1$.
\end{corollary}

Moreover, we also derive easily the $d$-dimensional extension of
\cite[Theorem 3]{BK}.

\begin{corollary}\label{nonconvex} {\bf (Boas -- Kac).}
Suppose that for some $n\in\NN$ the open set $\Omega\subset
\RR^d$ contains no integer multiples $kz$ of the point
$z\in\RR^d$ with $k>n$. Then we have again $ {\cal M}(\Omega,z)\le
M_n = 2 \cos{\pi \over n+2}$.
\end{corollary}

Apart from the convex case there are several cases of
\eqref{omegaz} when through the trigonometric extremal problem
\eqref{Mextr} either the precise value, or at least some estimate
can be found.

\begin{theorem}
\label{th:manycase}
Let $\Omega$ be a symmetric open set and $z\in\Omega$. Then the value of the
extremal quantity \eqref{omegaz} satisfies the following relations.

\begin{itemize}
\item[(i)]
If $H(\Omega,z)=\{n\}$, then ${\cal M}(\Omega,z)=\frac{1}{2\cos
{\pi \over 2n}}$.
\item[(ii)]
If $H(\Omega,z)=\NN_2\setminus\{n\}$, then
${\cal M}(\Omega,z)= \cos {\pi \over 2n}$.
\item[(iii)]
If $H(\Omega,z)=(n,\infty)\cap\NN_2$, then ${\cal M}(\Omega,z)=
\frac{1}{2\cos {\pi \over n+2}}$.
\item[(iv)]
If $H(\Omega,z)=2\NN+1$, then ${\cal M}(\Omega,z)= \frac{2}{\pi}$.
\item[(v)]
If $H(\Omega,z)=2\NN$, then ${\cal M}(\Omega,z)= \frac{\pi}{4}$.
\end{itemize}

\end{theorem}

\begin{remark} The extremal quantities $\cal M$ and $M$ are
monotonic in the sets $\Omega$ and $H$, respectively, hence the
above relations imply the corresponding inequalities when we know
only that e.g., $nz\in\Omega$, etc. We skip the formulation.
\end{remark}

\noindent {\bf Proof of Theorem \ref{th:manycase}.} In view of Theorem 
\ref{th:connection}, the calculation of 
$ {\cal M}(\Omega,z) $ hinges on finding the value of $M(H(\Omega,z))$. The
solutions of the corresponding trigonometric polynomial extremal problems,
relevant to the above list (i)-(v), can be looked up from the literutre as
follows. 
\begin{itemize}
\item[(i)]
An easy calculation, see e.g., \cite{R1}.
\item[(ii)]
See \cite{R1}, Proposition 1.
\item[(iii)]
See \cite{R3}.
\item[(iv)]
See the end of \cite{Sch}.
\item[(v)]
See \cite[p.\ 492-493]{R2}.
\end{itemize}

When ${\cal M}(\Omega,z)$ is known for a certain $H(\Omega,z)$, then further
cases can be obtained via the following duality result.

\begin{lemma}
\label{dualitylemma}
{\bf (see \cite{R1}).}
Let $H\subseteq\NN_2$ be arbitrary. Then we have
$$
M(H) M(\NN_2\setminus H)=2.
$$
\end{lemma}

In fact, this gives (ii) once (i) is known; (iii) and Corollary
\ref{nonconvex} and also
(iv) and (v) are similarly related, although they were obtained differently
in the works mentioned above.

To formulate the corresponding relation in Problem 1 we can record

\begin{corollary}
For any open set $\Omega\subseteq \RR^d$ and $z\in\Omega$ we have
$$
{\cal M}(\Omega,z){\cal M}(\Omega^*,z) = \frac 12 ,
$$
where $\Omega^*$ is any open, symmetric set containing $0$, $z$ and
$(\NN_2\setminus H(\Omega,z))z$, but disjoint from $ H(\Omega,z)z$.
\end{corollary}

Ending this section, let us recall that investigation of
Tur\'an-type problems started with keeping an eye on number
theoretic applications and connected problems. The interesting
paper of Gorbachev and Manoshina \cite{Gorbi} mentions \cite{KS}.

\begin{problem}\label{Deltaproblem} Determine
$$
\Delta (n) := \sup \{ M(H)/2 \, : \,\, H\subseteq \NN_2, |H|=n \}.
$$
\end{problem}

We only know (cf \cite{R1})
$$
1-{5 \over (n+1)^2} \le \Delta (n) \le 1 - {0.5 \over (n+1)^2}.
$$
The question is relevant to the Beurling theory of generalized primes, see
\cite{R4}.

\vskip0.6cm

{\bf \S 5 Proof of Theorem \ref{th:torus}}

\vskip0.6cm

As above, without loss of
generality we can restrict ourselves to sets $\Omega$ symmetric
about the origin. Similarly to the proof of Theorem
\ref{th:connection}, we are to prove two inequalities for both
cases.

\vskip0.4cm

\noindent \underline{Case $\#{\cal Z}=\infty \, :
\, {\cal M}^*(\Omega, z) \le M(H(\Omega,z))/2 $ }

\vskip0.4cm

Let $f\in {\cal F}^*(\Omega)\cap\ C^\infty(\TT^d)$. We consider
the measure
$$
\sigma_z^{(N)} := \sum_{k=-N}^{N} (1-\frac{|k|}{N})\delta_{kz}.
$$
This measure is positive definite since for all $n \in \ZZ^d$ we have
$$
\ft{\sigma_z^{(N)}}(n)=\int_{\TT^d} e^{-2\pi i \inner{n}{x}}
d\sigma_z^{(N)}(x)= \sum_{k=-N}^{N} (1-\frac{|k|}{N}) e^{2\pi i
\inner{n}{z}}=:K^{(N)}(2\pi \inner{n}{z}),
$$
where $K^{(N)}$ is the usual Fej\'er kernel, which is nonnegative.
Let us denote $H(N):=H(\Omega,z)\cap[2,N]$.

The function $f$ being continuous and even, the measure
\beql{rho}
\rho_z := f\cdot \sigma_z^{(N)} = f(0) \delta_0 +
\sum_{k \in \{1\}\cup H(N)}
(1-\frac{k}{N}) f(kz) (\delta_{kz}+\delta_{-kz})
\eeq
is well defined and, by $\ft{\rho_z}=\ft{f}*\ft{\sigma^{(N)}_z}$,
is positive definite as well.
In view of $f(0)=1$ we now find for arbitrary $n \in \ZZ^d$
that
$$
0 \le \widehat{\rho_z}(n) =
1+(2-{2\over N})f(z) \cos{2\pi\inner{z}{n}} +
  \sum_{k\in H(N)} (2-{2k\over N}) f(kz) \cos{2\pi k \inner{z}{n}}.
$$
Setting $t := \inner{z}{n}$ yields
$$
0 \le \varphi_N(t):=1 + 2(1-{1\over N})f(z)\cos{2\pi t} +
\sum_{k\in H(N)} 2 (1-{k\over N}) f(kz) \cos{2\pi k t}.
$$
Since $\#{\cal Z}=\infty $, here for the various values of $n
\in \ZZ^d$ the derived variable $t$ will be dense in $\TT$.

Hence we can conclude that in the infinite case $\varphi_N(t)\in
\Phi(H(\Omega,z))$. This gives $2(1-{1\over N})f(z) \le
M(H(\Omega,z))$ for all $N\in\NN$. Whence the stated inequality.

\vskip0.4cm

\noindent \underline{ Case $\#{\cal Z}=m<\infty \,:\, {\cal
M}^*(\Omega, z) \le M_m(H_m(\Omega,z))/2$ }

\vskip0.4cm

Let again $f\in {\cal F}^*(\Omega)\cap\ C^\infty(\TT^d)$. Now we
consider the measure
$$
\sigma_{z,m} := \frac 12 \sum_{k=-[{m-1\over 2}]}^{[{m-1\over 2}]}
\delta_{kz} + \frac 12 \sum_{k=-[{m\over 2}]}^{[{m\over 2}]}
\delta_{kz}.
$$
For all $n \in \ZZ^d$ we have
$$
\ft{\sigma_{z,m}}(n)=\int_{\TT^d} e^{-2\pi i \inner{n}{x}}
d\sigma_{z,m}(x)= 1 + \sum_{k=1}^{[{m-1\over 2}]} \cos {2\pi k
\inner{n}{z}} + \sum_{k=1}^{[{m\over 2}]} \cos {2\pi k
\inner{n}{z}}.
$$
Since $\#{\cal Z}=m < \infty $, where $m=[q_1,\dots,q_d]$ with
$z=({p_1 \over q_1},\dots,{p_d\over q_d}),
\,\, (p_j,q_j)=1 \,\, (j=1,\dots,d)$,
for the various values of $n \in
\ZZ^d$ the derived variable $t:=\inner {n}{z} $ will cover exactly
the values of $j/m \, ({\rm mod}\,\TT)$. For these values,
however, direct calculation shows that the above sum is either
exactly $m$ (in case $j\equiv 0\,({\rm mod}\,m)$), or vanishes.
Thus, again, the measure $\sigma_{z,m}$ will be positive definite.

The function $f$ being continuous and symmetric, the measure

\beql{rhom}
\rho_{z,m} := f\cdot \sigma_{z,m} = f(0) \delta_0 +
\sum_{k=1}^{[{m-1\over 2}]}  f(kz) (\delta_{kz}+\delta_{-kz}) +
\sum_{k=1}^{[{m\over 2}]} f(kz) (\delta_{kz}+\delta_{-kz})
\eeq

is well defined and, by
$\ft{\rho_{z,m}}=\ft{f}*\ft{\sigma_{z,m}}$, is positive definite
as well. In view of $f(0)=1$ we now find for all $n \in \ZZ^d$

\beql{rhohat} 0 \le \widehat{\rho_z}(n) = 1+2 f(z)\cos{2\pi t} +
\sum_{k=2}^{[{m-1\over 2}]} f(kz) \cos {2\pi k t} +
\sum_{k=2}^{[{m \over 2}]} f(kz) \cos {2 \pi k t}, \eeq
where $t=\inner{z}{n}$ as above. So let us write now
$$
\varphi_{z,m}(t):= 1+2 f(z)\cos{2\pi t} + \sum_{k=2}^{[{m-1\over
2}]} f(kz) \cos {2\pi k t} + \sum_{k=2}^{[{m\over 2}]} f(kz) \cos
{2\pi k t}.
$$
It follows that
$$
\varphi_{z,m}(t) = 1 + 2 f(z)\cos{2\pi t} + \sum_{k\in
H_m(\Omega,z)} c^*_k \cos{2\pi k t},
$$
for some $c^*_k \in \RR$.
Similarly as above, \eqref{rhohat} implies $\varphi_{z,m}(j/m)\ge
0 \, (j=0,\dots,m-1)$. That is, we conclude $\varphi_{z,m}\in \Phi
_m(H_m(\Omega,z))$ and thus $2 f(z) \le M_m(H_m(\Omega,z))$.
Hence the statement.

\vskip0.4cm

\noindent
\underline{ Case $\#{\cal Z}=\infty \, :
\, {\cal M}^*(\Omega,z) \ge  M(H(\Omega,z))/2$}

\vskip0.4cm

Let $\varphi$ be any trigonometric polynomial from the class
\eqref{Ficlass}. Then $\varphi$ has (restricted) spectral set $S$
and full spectrum  $S':=\{-1,0,1\}\cup \pm S$ with $S\subseteq
H:=H(\Omega,z)$ necessarily finite. Note that the supremum in the
definition \eqref{Mextr} of $M(H(\Omega, z))$ can be restricted
to the trigonometric polynomials of \eqref{Ficlass}.

Consider the measure
$$
\alpha_z = \delta_0 + (\lambda/2)(\delta_z + \delta_{-z}) +
\sum_{k\in S} (c_k/2) (\delta_{kz} + \delta_{-kz}),
$$
whose Fourier transform $\ft{\alpha_z}(n) =
\varphi(\inner{z}{n})\,(n\in\ZZ^d)$ is essentially the polynomial
$\varphi (t)$ itself. Hence $\alpha_z$ is a positive definite
measure.

Take now the ``triangle function'' $\Delta_{\epsilon}$, defined in
\eqref{triangle}, with a parameter $\epsilon$ small
enough to guarantee that
\begin{enumerate}
\item
The sets $kz + B_{\epsilon}$, $(k \in S')$, are disjoint, and
\item
These sets are all contained in $\Omega$, i.e., $\epsilon < {\rm
dist}\{\partial \Omega, S'z \}$.
\end{enumerate}
Since we consider only a finite subset $S$ of $H$, and
$S'=\{-1,0,1\} \cup \pm S)$, these conditions are met with some
positive $\epsilon$ as no two different multiples of $z$ are
equal in $\TT^d$. Finally define
$$
f := \alpha_z * \Delta_\epsilon,
$$
which is a positive definite function supported in $\Omega$ with
value $1$ at the origin and with $f(z) = \lambda / 2$. This
proves that ${\cal M}^*(\Omega, z) \ge \lambda /2$, hence taking
supremum over all polynomials $\varphi\in\Phi(H)$ concludes the
proof. \Qed

\vskip0.4cm

\noindent \underline{ Case $\#{\cal Z}=m<\infty \,:\, {\cal
M}^*(\Omega, z) \ge M_m(H_m(\Omega,z))/2  $ }

\vskip0.4cm

We denote here $H:=H_m(\Omega,z)$. Now take any $\varphi$ in
\eqref{Fimstar}.

Consider the measure
$$
\alpha_z = \delta_0 + (\lambda/2)(\delta_z + \delta_{-z}) +
\sum_{k<{m\over 2}, k \in H} (c_k/2) (\delta_{kz} +
\delta_{-kz})+c_{m/2} \delta_{mz/2},
$$
with the last term appearing only if $m$ is even and $m/2$
belongs to the spectral set \eqref{Hmstar}. Observe that for the
true spectrum of this measure we have

\beql{alphaspectrum}
S^*:=\supp \ft{\alpha_z}:=S^*(\alpha_z)
\subseteq \{-1,0,1\}\cup\pm H \setminus \{-m/2\}=
S'\setminus \{-m/2\}, \eeq
where the last term ($\setminus\{-m/2\}$) appears only if $m$ is even.
Thus it is easy
to see that the multiples $kz \,(k\in S^*)$ are different even in
$\TT^d$.

Now let us prove that $\alpha_z$ is positive definite. Taking
$n\in\ZZ^d$ arbitrarily, consider the Fourier transform
$$
\ft{\alpha_z}(n) = 1+\lambda \cos 2\pi\inner{z}{n} +
\sum_{k<{m\over 2}, k \in H} c_k \cos 2\pi k \inner{z}{n} +
c_{m/2} e^{-im\pi\inner{z}{n}}.
$$

Here, by the condition $\inner{z}{n}=j/m$ for some integer $j$, we have
in the last term $e^{-m\pi\inner{z}{n}}=(-1)^j=\cos\pi j=\cos
m\pi\inner{z}{n}$ and we get $\ft{\alpha_z}(n) = \varphi
(\inner{z}{n}) = \varphi(j/n)$. It follows that $\ft{\alpha_z}(n)
\ge 0$ by definition \eqref{Fimstar}.

Take now the ``triangle function'' $\Delta_{\epsilon}$ defined in
\eqref{triangle} with a parameter $\epsilon$ small enough to
ensure
\begin{enumerate}
\item
The sets $kz + B_{\epsilon}$, $(k \in S^*)$, are disjoint, and
\item
These sets are all contained in $\Omega$, i.e., $\epsilon < {\rm
dist} \{ \partial \Omega, S^*z \} $.
\end{enumerate}

These conditions are met with some positive $\epsilon$ since no
two different multiples $kz \, (k\in S^*)$ are equal in $\TT^d$,
and by definitions \eqref{Fimstar} and \eqref{alphaspectrum} we
necessarily have $S^*z\subseteq\Omega$.

Finally define
$$
f = \alpha_z * \Delta_\epsilon,
$$
which is a positive definite function supported in $\Omega$ with
value $1$ at the origin and with $f(z) = \lambda / 2$. This
proves that ${\cal M}^*(\Omega, z) \ge \lambda /2$, hence taking
supremum over all polynomials $\varphi\in\Phi_m(H)$ concludes the
proof. \Qed

\vskip0.6cm

{\bf \S 6 Applications of Theorem \ref{th:torus} and further connections}

\vskip0.6cm

Arestov, Berdysheva and Berens \cite{ABB} mention the one
dimensional symmetric interval special case of the following fact.

\begin{proposition}\label{lessthanstar}
Suppose $\Omega\subseteq (-\frac12,\frac12)^d$ is an open set. Then
$$
{\cal M}(\Omega,z)\le {\cal M}^*(\Omega,z) .
$$
\end{proposition}

{\bf Proof.} The original proof of \cite{ABB} uses the natural
periodization of functions $f\in {\cal F}(\Omega)$. Taking
$g(x):=\sum\limits_{n\in\ZZ^d}f(x-n)$ maps ${\cal F}(\Omega)$
injectively to ${\cal F}^*(\Omega)$, which proves the
Proposition. However, we have also an alternative argument here,
as Theorems \ref{th:connection} and \ref{th:torus} translate the
extremal problems in question to extremal problems for
trigonometric polynomials. In case \#${\cal Z}=\infty$ the $\RR^d$
and $\TT^d$ interpretations of \eqref{Homegaz} give
$H_{\RR^d}(\Omega,z)\subset
H_{\RR^d}(\Omega+\ZZ^d,z)=H_{\TT^d}(\Omega,z)$. For \#${\cal
Z}=m<\infty \quad H_{\RR^d}(\Omega,z)\subseteq [2,m-2]$. Indeed,
$-z\in \Omega\subseteq (-\frac12,\frac12)^d$, and as $0\ne mz$ but
$mz\equiv 0 \pmod{\TT^d}$, we obtain that $(m-1)z\notin\Omega$ in
$\RR^d$, and similarly for $k \ge m \quad kz\notin [-\frac 12,
\frac 12)^d$ excludes the possibility of $k \in
H_{\RR^d}(\Omega,z)$. Thus it is easy to see that

\beql{setcontroll} M_m(H_m(\Omega,z))=M_m(H(\Omega,z)\cap
[2,m-2]) =M_m(H_{\RR^d}(\Omega,z)) . \eeq

Now it is obvious that $\Phi_m(H)\supseteq\Phi(H)$ and thus
$M_m(H)\ge M(H)$ for arbitrary $H\subseteq\NN_2$, and we get the
assertion even for the finite case.

\begin{corollary}\label{upintegerpart}
Let $\Omega\subseteq (-\frac12,\frac12)^d$ be a convex, symmetric
domain. Then we have
$$
{\cal M}^*(\Omega,z) \ge w(||z||), \qquad {\rm where} \qquad
w(t):=\cos\frac{\pi}{\lceil 1/t \rceil +1}.
$$
\end{corollary}

{\bf Proof.} Corollary \ref{th:answer-to-pointwise-problem} gives
${\cal M}(\Omega,z) \ge w(||z||)$. Thus combining Proposition
\ref{lessthanstar} and Corollary
\ref{th:answer-to-pointwise-problem} proves the assertion.

\begin{remark}
The above estimate is a sharpening of (14) in \cite[Theorem 3]{ABB}.
\end{remark}

The following assertion is obvious both directly and by Theorem
\ref{th:connection}.

\begin{proposition}\label{dilateinvariance}
For all open sets $\Omega\subseteq \RR^d$ and $z\in\RR^d$, $\alpha>0$ we have
$$
{\cal M}(\alpha\Omega,\alpha z) = {\cal M}(\Omega,z) .
$$
\end{proposition}

\begin{proposition}\label{Ndilatation}
For  $\Omega\subseteq (-\frac12,\frac12)^d$ open, $z\in\TT^d$ and $N\in\NN$
we have
$$
{\cal M}^*(\frac1N\Omega,\frac1N z)\le {\cal M}^*(\Omega,z) .
$$
\end{proposition}

{\bf Proof.} One can work out the generalization of the proof of
\cite[Lemma 5]{ABB}, which is the one-dimensional interval special case
of this assertion. Instead, we note that $k\frac 1N z \in \frac
1N \Omega \pmod {\TT^d}$ entails $kz\in\Omega \pmod {\TT^d}$, and by
Theorem \ref{th:torus} the \#$\cal Z=\infty$ case follows.

On the other hand for finite \#${\cal Z}(z)=m<\infty$ we have
\#${\cal Z}(\frac 1N z) = Nm$ and
$\Phi_m(H)\supseteq\Phi_{mN}(H)$. Thus combining
\eqref{discreteformula} and \eqref{setcontroll} yields
$$
\begin{array} {rl}
2 {\cal M}^*(\Omega,z) = & M_m(H_m(\Omega,z))=
M_m(H_{\RR^d}(\Omega,z)
\\ = & M_m(H_{\RR^d}(\frac 1N \Omega, \frac 1N z) \ge
 M_{mN^*}(H_{\RR^d}(\frac 1N \Omega, \frac 1N z)
\\
= & M_{mN^*}(H_{mN^*}(\frac 1N \Omega, \frac1N z)) = 2{\cal
M}^*(\frac1N \Omega, \frac1N z).
\end{array}
$$
\Qed.

The next assertion is the generalization of
\cite[Theorem 4]{ABB}.

\begin{theorem}\label{limitcase}
For any bounded open set $\Omega\subset \RR^d$ and $z\in\RR^d$ we
have
$$
\lim_{\alpha \to +0}{\cal M}^*(\alpha\Omega,\alpha z) = {\cal
M}(\Omega,z).
$$
\end{theorem}

\begin{remark} Here the condition of boundedness ensures that
for $\alpha$ small enough we have $\alpha \Omega \subset
(-\frac12,\frac12)^d$ and the expression under the limit on the
left hand side is defined by \eqref{Mstar}.
\end{remark}

{\bf Proof.} Again, extending the original arguments of
\cite{Gorbi} or \cite{ABB} leads to a proof. There the idea is to
multiply $f\in {\cal F}^*(\alpha\Omega)$ by a fixed positive
kernel, say $\Delta_{\frac 14}$, and exploit that for $\alpha$
small $\Delta_{\frac 14} |_{\alpha\Omega}$ is approximately $1$.

Alternatively, we can argue as follows. Let
$\Omega$, be bounded by $R$, and let $\alpha<\frac{1}{2R}$: then $\alpha
\Omega \subseteq (-\frac12,\frac12)^d$. Moreover, using $\RR^d$
interpretation of the arising sets we always have
\beql{Hbound}
H_{\RR^d}(\Omega, z) = H_{\RR^d}(\alpha \Omega, \alpha z) \subset
\left[2,\frac {R}{|z|}\right],
\eeq
while $m(\alpha):=\#{\cal Z}(\alpha z) \ge \frac{1}{\alpha |z|}
\to \infty \quad (\alpha\to 0)$. Note that here for irrational
$\alpha$ we can have $m(\alpha)= + \infty$, but defining the
index function $m(\alpha)$ in this extended sense does not
question the asserted limit relation.

In what follows we unify terminology by writing
$H_{\infty}(\Theta,w)= H(\Theta, w)$ while keeping the notation
$H_n(\Theta,w)= H(\Theta, w)\cap [2,n/2]$ for finite $n$. For the
finite case we have $H_{\RR^d}(\alpha\Omega,\alpha z)=
H_{\RR^d}(\Omega, z)\subseteq [2,\frac{m(\alpha)}{2}]$, and in view
of \eqref{Homegaz} and \eqref{Hbound}
$H:=H_{m(\alpha)}(\alpha\Omega,\alpha z)=
H_{\TT^d}(\alpha\Omega,\alpha z)\cap [2,\frac{m(\alpha)}{2}]=
H_{\RR^d}(\alpha\Omega,\alpha z)=H_{\RR^d}(\Omega, z)$, too. Now
if $m(\alpha)=\infty$, then we are to consider the normalized,
nonnegative trigonometric polynomials $\varphi \in
\Phi_{\infty}(H):=\Phi(H) $ defined by \eqref{Ficlass}, while for
finite $m(\alpha)<\infty$, the function set to be considered is
$\Phi_m(H)$ defined by \eqref{Fimstar}.

Let now $\alpha_n\to 0$, and $\varphi_n$ be an extremal polynomial
in $\Phi_{m(\alpha_n)}(H)$. In view of the nonnegativity
conditions for these sets we get $|c_k|\le 2 \quad (k\in H)$,
applying finite Fourier Transform in case $m(\alpha_n)<\infty$.
Hence with $K:=\Ceil {\frac {2R}{|z|}}$ we find $\varphi_n \in
{\cal F}_K:= \{\varphi (t)=1+2\sum_{k=1}^K a_k \cos 2\pi kt
\given |a_k| \le 1, \, k=1,\dots ,K \}$, which is a compact
subset of $C(\TT)$. Thus without loss of generality we can
suppose that $\varphi_n \to \phi \in {\cal F}_K$ uniformly as $n
\to \infty$. Since $m(\alpha_n) \to \infty$, we must have $\phi
\ge 0$. Moreover, if we write $\phi (t)=1+2\sum_{k=1}^K a_k \cos
2\pi kt $ and $\varphi_n (t)=1+2\sum_{k=1}^K a_k^{(n)} \cos 2\pi
kt$, then $\lim_{n \to \infty} a_k^{(n)} = a_k$, so $\phi \in
\Phi(H)$ and
$$ \lim_{n \to \infty}{\cal M}^*(\alpha_n\Omega,\alpha_n z) =
\lim_{n \to \infty} a_1^{(n)} = a_1 \le {\cal M}(\Omega,z).$$

On the other hand Proposition \ref{lessthanstar} gives the
converse inequality. \Qed

\vskip0.6cm

{\bf \S 7 Calculations of extremal values for some special cases}

\vskip0.6cm

Now we formulate a periodic case analogon of the Boas-Kac result
Corollary \ref{BoasKac}.

\begin{proposition}\label{Zinomega}
Suppose that the open set  $\Omega\subseteq \TT^d$ contains all
integer multiples of the point $z\in\TT^d$, i.e., ${\cal Z} \subset
\Omega$ with ${\cal Z}$ defined in \eqref{Zset}. Then $ {\cal
M}^*(\Omega,z)=1$.
\end{proposition}

{\bf Proof.} In case \#$\Zset = \infty$, Theorem \ref{th:torus}
gives $\Msoz = M(H(\Omega,z))/2=M(\NN_2)/2=1$ immediately. Let now
\#$\Zset = m < \infty $. Then Theorem \ref{th:torus} yields $\Msoz
= M_m(H_m(\Omega,z))/2=M_m([2,m/2])/2$. To see that this quantity
achieves 1, it suffices to consider the cosine polynomial

$$
\varphi_m (t):=1+\sum_{k=1}^{[\frac{m-1}{2}]}\cos 2\pi kt +
\sum_{k=1}^{[\frac{m}{2}]}\cos 2\pi kt .
$$

Direct calculation proves again $\varphi_m (j/m) \ge 0 \quad
(j\in \NN)$, thus $\varphi_m \in \Phi_m([2,m/2])$ and now we find
$M_m([2,m/2])/2=1$. \Qed

With the following applications in mind we first prove

\begin{lemma}\label{evencase}
Let $m\in 2\NN$ be even. Then we have $M_m([2,m/2))=1+\cos
\frac{2\pi}{m}$.
\end{lemma}

{\bf Proof.} Let $m=2n$ and
$$
\varphi (t) = 1 + \sum_{k=1}^{n-1} c_k \cos 2\pi kt \in
\Phi_m(M_m([2,n))).
$$
Using the finite Fourier Transform coefficient formula and
$\varphi(j/m) \ge 0 \quad (j \in \NN)$ we obtain
\begin{eqnarray*}
c_1&=& \frac 2m \sum_{j=0}^{m-1} \varphi(\frac jm) \cos {\frac
{2\pi j}{m}} \\
& = & \frac 1n \sum_{l=0}^{n-1} \varphi(\frac ln) \cos {\frac {2\pi
l}{n}} + \frac 1n \sum_{l=0}^{n-1} \varphi(\frac {2l+1}{m}) \cos
(\frac {2\pi l}{n}+\frac {\pi} {n}) \\
& \le & \frac 1n \sum_{l=0}^{n-1} \varphi(\frac ln) + \frac 1n
\sum_{l=0}^{n-1} \varphi(\frac ln + \frac 1m) \cos (\frac
{\pi}{n}) = 1 + \cos (\frac {\pi}{n}).
\end{eqnarray*}

On the other hand take the cosine polynomial

$$
\phi_m (t):=1+\sum_{k=1}^{n-1} (1+ \cos \frac {\pi k} {n}) \cos
2\pi kt .
$$

Direct calculation gives

$$
\phi_m (\frac jm) = \left\{
\begin{array} {ll}
m & \qquad j\equiv 0 \pmod m \\
m/ 2 & \qquad j \equiv \pm 1 \pmod m \\
0 & \qquad {\rm otherwise},
\end{array} \right.
$$
whence $\phi_m (\frac jm) \ge 0 \quad (j \in \NN)$ and $\phi_m
\in \Phi_m(M_m(\left[2,n\right)))$. \Qed

\begin{corollary}\label{dimone}
{\bf (Arestov -- Berdysheva -- Berens \cite{ABB})} For dimension
one we have
\begin{itemize}
\item[(i)]
For $(p,q)=1$, $q$ even we have $\MM^*((-\frac 12, \frac12),\frac
pq) = \frac12 (1+\cos\frac{2\pi}{q})$.
\item[(ii)]
For $(p,q)=1$, $q$ odd we have $\MM^*((-\frac 12, \frac12),\frac
pq) = 1$.
\item[(iii)]
For $z \notin \QQ$ we have $\MM^*((-\frac 12, \frac12),z) = 1$.
\end{itemize}
\end{corollary}

{\bf Proof.} In case (i) \#$\Zset=q=2r$, and
$H(\Omega,z)=\NN_2\setminus r\NN$, $H_q^*(\Omega,z)=[2,r-1]$.
Hence in view of Theorem \ref{th:torus} it suffices to show that
$M_q^*(\left[2,r\right))=1+\cos (2\pi /q)$, which follows from
Lemma \ref{evencase}. For the cases (ii) and (iii) we clearly have
$\Zset \subseteq \Omega$, hence Proposition \ref{Zinomega} applies.
\Qed

Similarly to the above result of Arestov et al, we can also answer
the pointwise Tur\'an extremal problem for
$\Omega=(-\frac12,\frac12)^d$.

\begin{theorem}\label{intertorus}
Let $\Omega=(-\frac12,\frac12)^d \in \TT^d$. Then we have
\begin{itemize}
\item[(i)]
$\MM^*((-\frac 12, \frac12)^d,z) = 1$ if $z \notin \QQ^d$.
\end{itemize}
Moreover, if $z \in \QQ^d$,
$z=(\frac{p_1}{q_1},\dots,\frac{p_d}{q_d})$ with $(p_j,q_j)=1$,
$q_j=2^{s_j}t_j \quad (s_j\in\NN)$, $t_j\in 2\NN+1 \quad
(j=1,\dots,d)$
and $m:=[q_1,\dots,q_d]=2^st \quad t\in 2\NN+1$,
then we have either
\begin{itemize}
\item[(ii)]
$1\le s =s_1=\dots=s_d$, and then $\MM^*((-\frac 12, \frac12)^d,z)
= \frac 12 (1+\cos \frac{2\pi}{m})$, or
\item[(iii)]
$s=0$ or $\exists j$, $1\le j \le d$ with $s_j < s$ and then
$\MM^*((-\frac 12, \frac12)^d,z)=1$.
\end{itemize}
\end{theorem}

{\bf Proof.} Case (i) is covered by Proposition \ref{Zinomega}
above. If $z\in\QQ^d$, then the set defined in \eqref{Zset} is
finite and we have \#$\Zset=m=[q_1,\dots,q_d]$. Let us determine
the set $H(\Omega,z)$ first. For $k\in\NN$ we have $kz\notin
\Omega$ iff $kp_j/q_j\equiv 1/2 \pmod 1 \quad (j=1,\dots,d)$,
i.e., $2kp_j/q_j\equiv 1 \pmod 2 \quad (j=1,\dots,d)$. It follows
that $q_j|2k \quad (j=1,\dots,d)$, and we can not have a solution
$k\in\NN$ if $\exists j$ so that $q_j$ is odd, since then $2k/q_j$
must be even. Hence we can consider the case when all $s_j\ge1$
and, by $(p_j,q_j)=1$, all $p_j$ is odd. Then using $p_j\in
2\ZZ+1$ the condition becomes $2k/q_j\equiv 1 \pmod 2 \quad
(j=1,\dots,d)$. Hence $m=[q_1,\dots,q_d]|2k$ and $s=s_j\quad
(j=1,\dots,d)$ since otherwise for any $s_j<s$ we get
$2k/q_j=nm/q_j=n2^{s-s_j}t/t_j\equiv 0 \pmod 2$. In all, $kz
\notin \Omega$ occurs only in case (ii), while case (iii) will
again be covered by Proposition \ref{Zinomega}. In case (ii), when
$kz \notin \Omega$ happens, it occurs precisely for multiples of
$m/2 \in \NN$. That is, case (ii) now reduces to the determination
of $\Msoz=M_m(\left[2,m/2\right))/2=(1+\cos 2\pi/m )/2$ in view of
Theorem \ref{th:torus} and Lemma \ref{evencase}. \Qed


\noindent {\sc\small
Department of Mathematics, \\
University of Crete, \\
Knossos Ave., \\
714 09 Iraklio, Greece.\\
E-mail: {\tt mk@fourier.math.uoc.gr}\\
\ \\
and\\

\ \\
{\sc\small
Alfr\' ed R\' enyi Institute of Mathematics, \\
Hungarian Academy of Sciences, \\
1364 Budapest, Hungary}\\
E-mail: {\tt revesz@renyi.hu}
}

\end{document}